\documentclass[a4paper,oneside]{book}
\usepackage[english,russian]{babel}
\usepackage[cp1251]{inputenc}
\usepackage{graphicx}
\usepackage{amsmath}
\usepackage{amsfonts}
\usepackage{amssymb}
\usepackage{calc,ifthen}

\newcommand{\tezislarge}{\normalsize}
\newcommand{\tezissmall}{\small}
\providecommand{\No}{\textnumero}

\makeatletter
\renewcommand{\@makefnmark}{}
\newcommand{\l@author}{\@dottedtocline{2}{1.5em}{2.3em}}
\makeatother

\newcounter{first}

\newcommand{\tezis}[5][1]{\protect\setcounter{first}{\thepage}%
    \vspace{15mm plus 10mm minus 5mm}\begin{center}{\tezislarge\bfseries\uppercase{#2}}%
        \ifthenelse{\not\equal{#3}{}}{\footnote[1]{#3}}{}%
        \\*\textbf{#4}\addcontentsline{toc}{section}{\emph{#4,} #2}%
        \\*\textit{#5}%
    \end{center}\nopagebreak\setcounter{equation}{0}%
}

\newcommand{\athr}[2]{\footnote{#2}\addcontentsline{lof}{author}{#1}}

\newcommand{\abstrs}[6]{\vspace{-2ex}\tezissmall\ifthenelse{\not\equal{#1}{}}{%
        \ifthenelse{\not\equal{#3}{}}{УДК\ #3}{}
        \begin{quote}
            #1
            \ifthenelse{\not\equal{#2}{}}{\par\textit{Ключевые слова}:\ #2}{}%
        \end{quote}
    }{}%
    \selectlanguage{english}
    \begin{quote}
        \ifthenelse{\not\equal{#4}{}}{\begin{center}\textbf{#4}\end{center}}{}%
        #5
        \par\textit{Keywords}:\ #6
    \end{quote}

    \selectlanguage{russian}\tezislarge%
}

\begin{document}

\tezis{О равномерной устойчивости восстановления функций типа синуса с асимптотически отделенными нулями}
{Работа выполнена при финансовой поддержке РФФИ (проект \No\ 20-31-70005).}
{С.А. Бутерин}
{buterinsa@info.sgu.ru}

\athr{Бутерин С.А.}
{Бутерин Сергей Александрович, к.ф.-м.н., доцент,
 СГУ (Саратов, Россия);
 Sergey Buterin
 (Saratov State University, Saratov, Russia)}

\abstrs{Получена равномерная устойчивость восстановления целых функций специального вида по их нулям. К такому виду приводятся
характеристические определители усиленно регулярных дифференциальных операторов и пучков первого и второго порядков, включая
дифференциальные системы с асимптотически отделенными собственными значениями, характеристические числа которых лежат на прямой, содержащей
начало координат, и их нелокальных возмущений. Установлено, что зависимость таких функций от последовательностей их нулей является
липшицевой в естественных метриках на каждом шаре конечного радиуса. Результаты указанного типа могут быть использованы при исследовании
равномерной устойчивости обратных спектральных задач. Кроме того, получены общие теоремы об асимптотике нулей функций данного класса и об
их эквивалентном представлении в виде бесконечного произведения, которые дают соответствующие утверждения для многих конкретных операторов.
}
{функция типа синуса, усиленно регулярный дифференциальный оператор, собственные значения, характеристический определитель, равномерная
устойчивость, липшицева устойчивость}
{517.984}
{On the uniform stability of recovering sine-type functions with asymptotically separated zeros}
{We obtain a uniform stability of recovering entire functions of a special form from their zeros. To this form, one can reduce the
characteristic determinants of strongly regular differential operators and pencils of the first and the second orders, including
differential systems with asymptotically separated eigenvalues whose characteristic numbers lie on a line containing the origin, and their
non-local perturbations. We establish that the dependence of such functions on the sequences of their zeros possesses the Lipschitz
property with respect to natural metrics on each ball of a finite radius. Results of this type can be used for studying the uniform
stability of inverse spectral problems. In addition, general theorems on the asymptotics of zeros of functions of this class and on their
equivalent representation via an infinite product are obtained, which give the corresponding results for many specific operators.}
{sine-type function, strongly regular differential operator, eigenvalues, characteristic determinant, uniform stability, Lipschitz stability}

\vspace*{3mm}
{\bf\large 1. Введение}\\

В недавнее время появились результаты нового типа для классической обратной задачи Штурма--Лиувилля, относящиеся к вопросу ее равномерной
устойчивости~\cite{1}. Впоследствии в \cite{2, 3} была получена равномерная устойчивость обратных задач и для некоторых классов
интегро-дифференциальных операторов, причем использовался иной подход, нежели в \cite{1}, неотъемлемой частью которого стало доказательство
равномерной устойчивости восстановления характеристической функции рассматриваемого оператора по ее нулям, т.е. по спектру этого оператора.

Основной целью данной заметки является нахождение по возможности наиболее общего подкласса целых функций, на который последний результат
непосредственно обобщается, и которым охватываются характеристические определители достаточно широкого класса операторов, часто возникающих
в приложениях.

Установлено, что в качестве такого подкласса естественно рассмотреть всевозможные произведения алгебраических полиномов и функций типа
синуса, нули которых асимптотически отделены.

Помимо равномерной устойчивости восстановления таких функций по их нулям, в работе получены общие теоремы об асимптотике нулей и об
эквивалентном представлении функций этого класса в виде бесконечного произведения, которые дают соответствующие утверждения для многих
конкретных операторов без необходимости воспроизводить полное доказательство в каждом частном случае. Полученные результаты применимы для
исследования различных аспектов прямых и обратных спектральных задач и, в том числе, равномерной устойчивости последних.

Для начала проиллюстрируем равномерную устойчивость характеристической функции на примере краевой задачи Штурма--Лиувилля
\begin{equation}\label{1.1}
-y''+q(x)y=\lambda y, \quad 0<x<\pi, \quad y(0)=y(\pi)=0,
\end{equation}
где  $q(x)$ -- комплекснозначная функция из $L_2(0,\pi).$ Собственные значения задачи (\ref{1.1}) с учетом кратности совпадают с нулями
$\lambda_n,$ $n\in{\mathbb N},$ ее характеристической функции $\Delta(\lambda)=Y(\pi,\lambda),$ где $Y(x,\lambda)$ -- решение уравнения в
(\ref{1.1}), удовлетворяющее начальным условиям $Y(0,\lambda)=0$ и $Y'(0,\lambda)=1.$ Без ущерба для общности считаем, что среднее значение
$q(x)$ на интервале $(0,\pi)$ обращается в ноль. Тогда
\begin{equation}\label{1.2}
\Delta(\lambda)=\frac{\sin z\pi}{ z}+ \int_0^\pi u(x)\frac{\cos z x}{ z^2}\,dx, \quad u(x)\in L_2(0,\pi), \; \int_0^\pi u(x)\,dx=0, \;\;
z^2=\lambda,
\end{equation}
откуда, в частности, вытекает известная асимптотика $\{\lambda_n-n^2\}\in l_2.$ Наряду с $\Delta(\lambda)$ рассмотрим еще одну функцию
$\tilde\Delta(\lambda)$ вида (\ref{1.2}) с нулями $\tilde\lambda_n,$ для которой символ $\tilde\omega$ будет иметь тот же смысл, что и
$\omega$ -- для $\Delta(\lambda).$ Одним и тем же символом $C_r$ будем обозначать {\it различные} положительные константы, зависящие только
от~$r.$ Следующая теорема дает равномерную устойчивость восстановления $u(x)$~по нулям~$\Delta(\lambda).$

\medskip
{\bf Теорема 1. \cite{3} }{\it Для всякого $r>0$ справедлива оценка
$$
\|u-\tilde u\|_{L_2(0,\pi)}\le C_r\|\{\lambda_n-\tilde\lambda_n\}\|_{l_2},
$$
коль скоро $\|\{\lambda_n-n^2\}\|_{l_2}\le r$ и $\|\{\tilde\lambda_n-n^2\}\|_{l_2}\le r.$ }

\medskip
Аналогичное утверждение имеет место и в случае, когда потенциал $q(x)$ является комплексной функцией-распределением из
$W_2^{-1}[0,\pi].$ Тогда характеристическая функция, как известно, принимает вид (см., например, \cite{4})
\begin{equation}\label{1.3}
\Delta(\lambda)=\frac{\sin z\pi}{ z}+ \int_0^\pi v(x)\frac{\sin z x} z\,dx, \quad v(x)\in L_2(0,\pi),
\end{equation}
и, как следствие, $\{\lambda_n/n-n\}\in l_2.$

\medskip
{\bf Теорема 2. \cite{2} }{\it Для всякого $r>0$ справедлива оценка
\begin{equation}\label{1.3-0}
\|v-\tilde v\|_{L_2(0,\pi)}\le C_r\Lambda_1(\{\lambda_n\},\{\tilde\lambda_n\}), \quad \Lambda_j(\{\lambda_n\},\{\tilde\lambda_n\})
:=\Big\|\Big\{\frac{\lambda_n-\tilde\lambda_n}{n^j}\Big\}\Big\|_{l_2},
\end{equation}
коль скоро $\Lambda_1(\{\lambda_n\},\{n^2\})\le r$ и $\Lambda_1(\{\tilde\lambda_n\},\{n^2\})\le r.$ }

\medskip
В \cite{2} получена равномерная устойчивость $v(x)$ и в равномерной норме. А именно, в условиях теоремы~2 (с той же метрикой
$\Lambda_1(\,\cdot\,,\,\cdot\,))$ справедлива также оценка
$$
\|v-\tilde v\|_{L_\infty(0,\pi)}\le C_r\Big\|\Big\{\frac{\lambda_n-\tilde\lambda_n}{n}\Big\}\Big\|_{l_1},
$$
обе части которой могут обращаться в бесконечность.

Заметим, что речь идет о равномерной устойчивости восстановления не только функций $u(x)$ и $v(x),$ но и самой $\Delta(\lambda).$ В самом
деле, применяя равенство Парсеваля для преобразования Фурье, можно переформулировать (и объединить) теоремы~1 и~2 в следующем эквивалентном
виде.

\medskip
{\bf Теорема 3. }{\it Пусть $j\in\{0,1\}.$ Тогда для любого $r>0$ справедлива оценка
$$
\|\theta_j\|_{L_2(-\infty,\infty)}\le C_r\Lambda_j(\{\lambda_n\},\{\tilde\lambda_n\}), \quad
\theta_j(x):=x^{2-j}(\Delta-\tilde\Delta)(x^2),
$$
коль скоро $\Lambda_j(\{\lambda_n\},\{n^2\})\le r$ и $\Lambda_j(\{\tilde\lambda_n\},\{n^2\})\le r,$ где $\Lambda_j(\,\cdot\,,\,\cdot\,)$
определено в (\ref{1.3-0}). }

\medskip
Также отметим, что $f(z)= z\Delta( z^2)$ является функцией типа синуса  (см., например, \cite{5-1}), а ее нули асимптотически просты и
отделены. Естественно предположить, что результаты аналогичные теоремам~1--3 можно получать и для других классов операторов,
характеристические функции которых в комплексной плоскости надлежащей степени спектрального параметра обладают аналогичными свойствами.

В то время как доказательства теорем~1 и~2 аналогичны, никакая из них не следует непосредственно из другой. Мотивацией к данной работе
первоначально послужила потребность унифицировать подход к доказательству и по возможности распространить его на более общий класс функций.

В следующем разделе дается описание искомого класса и формулируются основные результаты работы. Доказательства приведены в разделах~3--5.
\\

{\bf\large 2. Основные результаты}
\\

Зафиксируем $N\in{\mathbb N}\cup\{0\},$ а также $b>0$ и рассмотрим целую функцию вида
\begin{equation}\label{2.1}
\theta(z)=S(z)+\int_{-b}^b w(x)\exp(izx)\,dx, \quad S(z)=P_N(z)S_0(z), \quad w(x)\in L_2(-b,b),
\end{equation}
где $P_N(z)$ является алгебраическим многочленом степени $N,$ а $S_0(z)$ -- функцией типа синуса экспоненциального типа $b,$ нули которой
$z^0_n,$ $n\in{\mathbb N},$ асимптотически просты и отделены, т.е. $\displaystyle\inf|z^0_n-z^0_k|>0$ при $n\ne k $ и $n,k\gg1.$

Напомним, что целая функция $S_0(z)$ экспоненциального типа называется функцией типа синуса (типа $b),$ если найдутся положительные
константы $c,$ $C$ и $K,$ при которых выполняется двусторонняя оценка
\begin{equation}\label{2.2}
c<|S_0(z)|\exp(-|{\rm Im\,}z|b)<C, \quad |{\rm Im\,}z|>K.
\end{equation}
Заметим, что корректировкой $P_N(z)$ и $S_0(z)$ можно избавиться от интегрального слагаемого в (\ref{2.1}), однако его выделение необходимо
для наших целей. Отметим также, что к виду $P_N(z)S_0(z)$ может быть приведена и всякая целая функция вида
$$
S(z)=\sum_{j=0}^N z^{N-j}s_j(z), \quad s_j(z)=O(\exp(|{\rm Im\,}z|b), \quad z\to\infty, \quad j=\overline{0,N},
$$
где $s_0(z)$ -- некоторая функция типа синуса с асимптотически отделенными нулями.

Легко видеть, что функции (\ref{1.2}) и (\ref{1.3}) после умножения на соответствующую степень $z$ принимают вид (\ref{2.1}) с $b=\pi$ при
$N=1$ и $N=0,$ соответственно. Аналогичным образом к виду (\ref{2.1}) могут быть приведены характеристические определители усиленно
регулярных \cite{6} дифференциальных операторов и пучков первого и второго порядков на конечном интервале, включая дифференциальные системы
с асимптотически отделенными собственными значениями, характеристические числа главной части которых располагаются на одной прямой,
содержащей и начало координат, и их всевозможных нелокальных возмущений.

Отметим, что требование {\it усиленной} регулярности существенно. Например, оператор Штурма--Лиувил\-ля с периодическими либо с
антипериодическими краевыми условиями, будучи регулярным, но не усиленно регулярным и обладающим асимптотически двукратными собственными
значениями, сюда не входит и требует отдельного исследования. Данное ограничение связано со следующим обстоятельством. В силу известной
теоремы Б.\,Я.~Левина \cite{4-1}, система экспонент $\{\exp(iz^0_nx)\}_{n\ge1}$ является базисом в $L_2(-b,b),$ если все нули
$\{z^0_n\}_{n\ge1}$ отделены. Впоследствии В.\,Д.~Головиным \cite{4-2} было уточнено, что соответствующий базис является базисом Рисса (см.
тж. \cite{5-1,5}). Эти утверждения легко обобщаются и для асимптотически отделенных нулей путем той или иной коррекции конечного числа
экспонент, соответствующих кратным нулям, как сделано, например, ниже в леммах~1 и~2. Однако случай асимптотически неотделенных и, в
частности, бесконечного числа кратных нулей требует отдельного рассмотрения.

Отметим, что верхняя оценка в (\ref{2.2}) необходимо верна и для всех $z\in{\mathbb C};$ см., например, \cite{5}. Там же показано, что
внутрь полосы может быть продолжена и нижняя оценка, но следующим образом:
\begin{equation}\label{2.2-1}
|S_0(z)|\ge c_\delta\exp(|{\rm Im\,}z|b), \quad {\rm dist}(z,\{z^0_n\}_{n\ge1})\ge\delta,
\end{equation}
где $c_\delta>0$ зависит от $\delta>0.$ Простым следствием этих уточнений и принципа максимума для аналитических функций является следующая
двусторонняя оценка, также приведенная в \cite{5}. А именно, найдутся такие $N_1$ и $N_2,$ что
\begin{equation}\label{2.3}
0<N_1<|S_0'(z^0_n)|<N_2<\infty,
\end{equation}
коль скоро $z^0_n$ является простым нулем функции $S_0(z).$ Введем обозначение
\begin{equation}\label{2.3-1}
\mu_n:=\left\{\begin{array}{cl}z_n^0, & z_n^0\ne0,\\
-1, & z_n^0=0,
\end{array}\right. \quad n\ge1-N,
\end{equation}
где $\{z_n^0\}_{n=\overline{1-N,0}}$ -- нули полинома $P_N(z).$

В работе получены следующие теоремы.

\medskip
{\bf Теорема 4. }{\it Всякая функция $\theta(z)$ вида (\ref{2.1}) обладает бесконечным множест\-вом нулей $\{z_n\}_{n\ge1-N},$ которые, в
свою очередь, имеют вид
\begin{equation}\label{2.4}
z_n=z^0_n+\frac{\varkappa_n}{\mu_n^N}, \quad \{\varkappa_n\}\in l_2.
\end{equation}
}

\medskip
{\bf Теорема 5. }{\it При фиксированной $S(z),$ функция $\theta(z)$ однозначно определяется заданием всех своих нулей за исключением любых
$N$ штук, т.е. заданием последовательности $\{z_n\}_{n\ge1}.$

Кроме того, справедливо представление
\begin{equation}\label{2.6}
\theta(z)=\alpha\exp(\beta z)\prod_{n=1-N}^\infty\frac{z_n-z}{\mu_n}\exp\Big(\frac{z}{\mu_n}\Big),
\end{equation}
где $\beta=s+\gamma,$ а $s$ -- кратность нуля $S(z)$ в точке ноль, и
\begin{equation}\label{2.5}
\alpha=\lim_{z\to0}\frac{S(z)}{z^s}, \quad \gamma=\lim_{z\to0}\frac{d}{dz}\ln\frac{S(z)}{z^s}.
\end{equation}
}

\medskip
Следующая теорема дает обратное утверждение к теоремам~4 и~5.

\medskip
{\bf Теорема 6. }{\it Зафиксируем $S(z)$ в (\ref{2.1}). Тогда для любой комплексной последовательности $\{z_n\}_{n\ge1}$ вида (\ref{2.4})
существует единственный набор чисел $\{z_n\}_{n=\overline{1-N,0}},$ такой что функция $\theta(z),$ определенная формулами (\ref{2.6}) и
(\ref{2.5}), имеет вид (\ref{2.1}).

Если числа $\{z_n\}_{n=\overline{1-N,0}}$ также выбраны произвольным образом, то соответствующая функция $\theta(z)$ примет вид
\begin{equation}\label{2.6-1}
\theta(z)=S(z)+P_{N-1}(z)S_0(z)+\int_{-b}^b \tilde w(x)\exp(izx)\,dx, \quad \tilde w(x)\in L_2(-b,b),
\end{equation}
где $P_{N-1}(z)$ -- некоторый многочлен степени меньшей $N.$ }

\medskip
Сформулируем теперь главный результат работы, для чего наряду с $\theta(z)$ будем рассматривать еще одну функцию $\tilde\theta(z)$
того же вида (\ref{2.1}) и с той же главной частью $S(z),$ но с другой подынтегральной функцией $w(x):$
$$
\tilde\theta(z)=S(z)+\int_{-b}^b \tilde w(x)\exp(izx)\,dx, \quad \tilde w(x)\in L_2(-b,b).
$$
Условимся, что если некоторый символ $\omega$ обозначает объект, относящийся к $\theta(z),$ то тот же символ с тильдой $\tilde\omega,$
будет обозначать аналогичный объект, соответствующий функции $\tilde\theta(z),$ и $\hat\omega:=\omega-\tilde\omega.$ Имеет место следующая
теорема.

\medskip
{\bf Теорема 7. }{\it Для любого $r>0$ справедлива оценка
\begin{equation}\label{2.8}
\|\hat w\|_{L_2(-b,b)}\le C_r\|\{\mu_n^N\hat z_n\}_{n\ge1-N}\|_{l_2},
\end{equation}
коль скоро $\|\{\mu_n^N(z_n-z^0_n)\}_{n\ge1-N}\|_{l_2}\le r$ и $\|\{\mu_n^N(\tilde z_n-z^0_n)\}_{n\ge1-N}\|_{l_2}\le r.$ }

\medskip
{\bf Замечание 1.} Поскольку, согласно теореме~5, при восстановлении $w(x)$ любые $N$ нулей функции $\theta(z)$ являются лишними, есть
основание ожидать, что некоторый аналог теоремы~7 будет иметь место и после замены “$n\ge1-N$'' на “$n\ge1$''. Однако получение строгой
формулировки выходит за рамки настоящей работы.

\medskip
Аналогично теореме~3 неравенство (\ref{2.8}) при помощи равенства Парсеваля для преобразования Фурье приводится к эквивалентному виду:
\begin{equation}\label{2.10}
\|\hat\theta\|_{L_2(-\infty,\infty)}\le C_r\|\{\mu_n^N\hat z_n\}_{n\ge1-N}\|_{l_2},
\end{equation}
что означает равномерную устойчивость восстановления непосредственно~$\theta(z).$

Кроме того, в силу теоремы Планшереля и Полиа (см. \cite{5}, с.~50), имеем
$$
\int_{-\infty}^\infty |f(x+iy)|^p\,dx\le\|f\|_{L_p(-\infty,\infty)}^pe^{pb|y|}, \quad p>0,
$$
для всякой функции $f(z)$ экспоненциального типа $b$, $p$-я степень которой суммируема на вещественной оси. Таким образом, (\ref{2.10})
влечет оценку
$$
\|\hat\theta(\,\cdot\,+iy)\|_{L_2(-\infty,\infty)}\le C_r\|\{\mu_n^N\hat z_n\}_{n\ge1-N}\|_{l_2},
$$
равномерную по $y$ на ограниченных множествах.

Несмотря на то, что функция $\theta(z)$ вида (\ref{2.1}) является функцией типа синуса тогда и только тогда, когда $N=0,$ можно было бы в
полной мере оправдать название статьи, ограничившись только частным случаем, когда $z_k=\tilde z_k$ фиксированы, например, при
$k=\overline{1-N,0},$ а вместо функций $\theta(z)$ и $\tilde\theta(z)$ рассмотреть функции
$$
\theta_1(z):=\theta(z)\prod_{n=1-N}^0\frac1{z-z_k}, \quad \tilde\theta_1(z):=\tilde\theta(z)\prod_{n=1-N}^0\frac1{z-z_k}.
$$
Тогда в теореме~7 все вхождения ''$n\ge1-N$'' автоматически заменятся на ''$n\ge1$''. Однако это ограничило бы общность, поскольку в
соответствующих представлениях характеристических функций различных операторов часто присутствует степенной множитель -- например, в случае
оператора Штурма--Лиувилля с краевыми условиями Неймана. Кроме того, иногда удобнее наоборот -- добавлять фиктивные нули для конкретных
характеристических функций, чтобы избавляться от $z$ в знаменателе.
\\

{\bf\large 3. Доказательство теоремы 4}
\\

Прежде чем перейти непосредственно к доказательству, проведем некоторую подготовительную работу. Будем иметь
\begin{equation}\label{3.1}
f(z):=\theta(z)-S(z)=\int_{-b}^b w(x)\exp(izx)\,dx =o(\exp(|{\rm Im\,}z|b)), \quad z\to\infty.
\end{equation}

Нам потребуется следующее вспомогательное утверждение.

\medskip
{\bf Лемма 1. }{\it Для всякой ограниченной последовательности $\{\alpha_n\}_{n\ge1}$ имеет место включение
$\{f(z^0_n+\alpha_n)\}_{n\ge1}\in l_2,$ где $\{z^0_n\}_{n\ge1}$ -- последовательность нулей $S_0(z).$ }

\medskip
{\it Доказательство.} Согласно (\ref{3.1}) справедливо соотношение
$$
f(z^0_n+\alpha_n)=\sum_{\nu=0}^\infty\frac{(i\alpha_n)^\nu}{\nu!}\int_{-b}^b x^\nu w(x)\exp(iz^0_nx)\,dx,
$$
откуда при помощи неравенства Коши--Буняковского получаем оценку
\begin{equation}\label{3.2}
|f(z^0_n+\alpha_n)|^2\le C\sum_{\nu=0}^\infty\frac1{\nu!}\Big|\int_{-b}^b x^\nu w(x)\exp(iz^0_nx)\,dx\Big|^2, \quad
C:=\sup_{n\ge1}e^{|\alpha_n|^2}<\infty.
\end{equation}
Как уже отмечалось в предыдущем разделе, поскольку $S_0(z)$ является функцией типа синуса, а ее нули $\{z^0_n\}_{n\ge1}$ асимптотически
отделены, всякая система функций $\{\exp(iz^1_nx)\}_{n\ge1}$ образует базис Рисса в $L_2(-b,b),$ если последовательность
$\{z^1_n\}_{n\ge1}$ получена из $\{z^0_n\}_{n\ge1}$ после любой замены конечного числа кратных элементов, чтобы только
$\{\exp(iz^1_nx)\}_{n\ge1}$ была полна в $L_2(-b,b).$ Следовательно, существуют положительные константы $m$ и $M,$ такие что
\begin{equation}\label{3.3}
m\|g\|_{L_2(-b,b)}^2 \le \sum_{n=1}^\infty\Big|\int_{-b}^b g(x)\exp(iz^1_nx)\,dx\Big|^2 \le M\|g\|_{L_2(-b,b)}^2
\end{equation}
для любой функции $g(x)\in L_2(-b,b).$ Применяя неравенства (\ref{3.2}) и (\ref{3.3}), получаем
$$
\sum_{n=k}^\infty|f(z^0_n+\alpha_n)|^2\le C\sum_{\nu=0}^\infty\frac1{\nu!}\sum_{n=1}^\infty\Big|\int_{-b}^b x^\nu
w(x)\exp(iz^1_nx)\,dx\Big|^2 \qquad\qquad\qquad
$$
$$
\qquad\qquad\qquad\qquad\qquad\qquad\qquad\qquad\qquad\le CMe^{b^2}\|w\|_{L_2(-b,b)}^2<\infty,
$$
где $k$ выбрано таким образом, чтобы исключить $z^0_n\ne z^1_n.$ Лемма доказана.

\medskip
Для доказательства теоремы заметим, что в силу асимптотической отделимости нулей $S(z)$ существует $\delta>0,$ при котором найдется система
неограниченно расширяющихся замкнутых контуров $\Gamma_n,$ $n\ge1,$ в плоскости переменной $z,$ такая что ${\rm dist}(\{\Gamma_n\}_{n\ge1},
\{z^0_n\}_{n\ge1-N})\ge\delta.$ Обозначим через $k_n$ число нулей функции $S(z)$ внутри контура $\Gamma_n:$ $z^0_{1-N},\ldots,z^0_{k_n-N}.$
Тогда, в силу оценок (\ref{2.2-1}) и (\ref{3.1}), а также теоремы Руше, ${\rm int\,}\Gamma_n$ при достаточно больших $n$ содержит в
точности $k_n$ нулей функции $\theta(z):$ $z_{1-N},\ldots,z_{k_n-N}.$ При этом, согласно (\ref{3.1}), можно выбрать такое  $n=n_\delta,$
чтобы всюду в ${\rm ext\,}\Gamma_n$ выполнялось строгое неравенство $|f(z)|<c_\delta|P_N(z)|\exp(|{\rm Im\,}z|b),$ где $c_\delta$ взято из
оценки (\ref{2.2-1}). Тогда в каждом круге $\gamma_\delta(z^0_k):=\{z:|z-z^0_k|<\delta\},$ лежащем в ${\rm ext\,}\Gamma_n,$ содержится в
точности один нуль функции $\theta(z).$ Никаких других нулей у функции $\theta(z)$ нет, в чем можно легко убедиться, предположив противное
и рассмотрев контур $\Gamma_{n_1}$ с достаточно большим $n_1>n.$ Таким образом, приходим к бесконечной последовательности
$\{z_n\}_{n\ge1-N}$ нулей функции $\theta(z),$ которые в силу произвольной малости $\delta>0$ имеют вид $z_n=z^0_n+\varepsilon_n,$ где
$\varepsilon_n=o(1)$ при $n\to\infty.$

Далее, в силу (\ref{3.1}) имеет место соотношение $S(z_n)=-f(z_n),$ $n\ge1-N,$ откуда согласно лемме~1 имеем $\{S(z_n)\}\in l_2.$ С другой
стороны, $S(z_n)=\varepsilon_nS_n(z_n),$ где функция $S_n(z):=(z-z^0_n)^{-1}S(z)$ доопределена до целой. Таким образом, в силу (\ref{2.1})
и (\ref{2.2-1}) для достаточно больших $n$ будем иметь
$$
|S_n(z_n)|>\min_{|z-z^0_n|=\delta}\frac{|S(z)|}{|z-z^0_n|}\ge\frac{c_\delta}\delta\min_{|z-z^0_n|=\delta}|P_N(z)|\ge C_\delta|\mu_n|^N>0,
$$
что вместе с предыдущими рассуждениями дает $\{(\mu_n)^N\varepsilon_n\}\in l_2.$ Теорема доказана.
\\

{\large\bf Доказательство теорем~5 и~6}
\\

Сначала докажем следующее вспомогательное утверждение.

\medskip
{\bf Лемма 2. }{\it Пусть $S_0(z)$ -- некоторая функция типа синуса типа $b$ с асимптотически отделенными нулями $\{z^0_n\}_{n\ge1},$ a
$\{\kappa_n\}_{n\ge1}$ -- произвольная последовательность из $l_2.$ Обозначим $z_n:=z^0_n+\kappa_n,$ $n\in{\mathbb N},$ считая для
удобства, что кратные $z_n$ занумерованы подряд: $z_n=\ldots=z_{n+m_n-1},$ где $m_n$ -- кратность $z_n$ в последовательности
$\{z_k\}_{k\ge1}.$ Положим
$$
\sigma:=\{n:n\in{\mathbb N},z_n\ne z_{n-1}, n\ge2\}\cup\{1\},
$$
$$
e_{k+\nu}(x):=x^\nu\exp(iz_kx),\;\; k\in\sigma, \;\;\nu=\overline{0,m_k-1}.
$$
Тогда система функций $\{e_n(z)\}_{n\ge1}$ образует базис Рисса в $L_2(-b,b).$ }

\medskip
{\it Доказательство.} Нетрудно видеть, что система $\{e_n(z)\}_{n\ge1}$ квадратично близка к некоторому базису Рисса --  например, к базису
$\{\exp(iz^1_nx)\}_{n\ge1}$ в доказательстве леммы~1. Остается заметить, что полнота $\{e_n(z)\}_{n\ge1}$ вытекает из оценки
$$
\frac1{S_0(z)}\prod_{n=1}^\infty\frac{z-z_n^0}{z-z_n}\int_{-b}^b g(x)\exp(izx)\,dx =o(1),\qquad\qquad\qquad\qquad\qquad\qquad
$$
\vspace*{-7mm}
$$
\qquad\qquad\qquad\qquad\qquad\qquad\qquad\qquad z\to\infty, \quad {\rm dist}(z,\{z_n\}_{n\ge1})\ge\delta>0,
$$
выполняющейся, в частности, для любой функции $g(x)\in L_2(-b,b)$ и являющейся очевидным следствием неравенства (\ref{2.2-1}) вместе с
оценкой
$$
G(z)=O(1), \quad z\to\infty, \quad {\rm dist}(z,\{z_n\}_{n\ge1})\ge\delta, \quad G(z):=\prod_{n=1}^\infty\frac{z-z_n^0}{z-z_n}.
$$
Раз для конечного произведения последняя оценка очевидна, для ее обоснования в общем случае достаточно считать $|\kappa_n|\le\delta/2$ при
всех $n\in {\mathbb N}.$ Тогда можно записать
$$
G(z)=\exp\Big(\sum_{n=1}^\infty \ln\Big(1+\frac{\kappa_n}{z-z_n}\Big)\Big) =\exp\Big(\sum_{n=1}^\infty \sum_{\nu=0}^\infty
\frac{(-1)^\nu}{\nu+1}\Big(\frac{\kappa_n}{z-z_n}\Big)^{\nu+1}\Big),
$$
откуда получаем $|G(z)|\le\exp(2H(z)),$ где, в силу неравенства Коши--Буняковского,
$$
H(z):=\sum_{n=1}^\infty\frac{|\kappa_n|}{|z-z_n|} \le\|\{\kappa_n\}\|_{l_2}\sqrt{\sum_{n=1}^\infty\frac1{|z-z_n|^2}}.
$$
Поскольку в произвольной вертикальной полосе фиксированной ширины число точек~$z_n$ равномерно ограничено, для доказательства
ограниченности последней суммы при ${\rm dist}(z,\{z_n\}_{n\ge1})\ge\delta$ достаточно установить оценку
$$
H_1(z):=\sum_{k=-\infty}^\infty\frac1{|z-k-\alpha_k|^2}\le C_\delta, \quad {\rm dist}(z,\{k+\alpha_k\}_{k\in{\mathbb Z}})\ge\delta,
$$
где последовательность $\{\alpha_k\}_{k\in{\mathbb Z}}$ ограничена по модулю некоторой константой~$M.$ Пусть ${\rm Re\,}\,z\in(-1/2,1/2],$
что всегда можно достичь сдвигом индекса $k\to k+k_0$ и заменой $\alpha_k\to\alpha_{k+k_0}.$  Выберем $n_0$ так, чтобы $|z-k|\ge 2M$ при
$|k|\ge n_0.$ Тогда
$$
H_1(z)\le\frac{2n_0-1}{\delta^2} +\sum_{|k|\ge n_0}\frac1{(|z-k|-M)^2} \le\frac{2n_0-1}{\delta^2} +\sum_{|k|\ge n_0}\frac4{|z-k|^2}\le
C_\delta,
$$
причем $n_0,$ очевидно, не зависит от $k_0.$ Лемма доказана.

\medskip
Приведем доказательство теоремы~5. В силу леммы~2 функция $w(x)$ в представлении (\ref{2.1}) однозначно определяется соотношениями
\begin{equation}\label{4.2-1}
-S^{(\nu)}(z_k)=\int_{-b}^b  w(x)(ix)^\nu\exp(iz_kx)\,dx, \quad k\in\sigma, \quad \nu=\overline{0,m_k-1},
\end{equation}
не содержащими $z_n$ при $n<1,$ откуда следует первая часть утверждения теоремы.

Далее, согласно теореме Адамара о разложении целой функции конечного порядка в бесконечное произведение, будем иметь
\begin{equation}\label{4.1}
\theta(z)=C_0z^p\exp(C_1z)\prod_{z_n\ne0}\Big(1-\frac{z}{z_n}\Big)\exp\Big(\frac{z}{z_n}\Big), \quad
\end{equation}
где $p$ -- кратность нуля функции $\theta(z)$ в нуле. В частности, имеем
\begin{equation}\label{4.2}
S(z)=\alpha z^s\exp(\gamma z)\prod_{z^0_n\ne0}\Big(1-\frac{z}{z^0_n}\Big)\exp\Big(\frac{z}{z^0_n}\Big),
\end{equation}
где числа $\alpha$ и $\gamma$ определены в (\ref{2.5}). Вычислим $C_0$ и $C_1$ в (\ref{4.1}). В силу (\ref{2.2-1}) и~(\ref{3.1}) имеем
$\theta(z)/S(z)\to1$ при $z^2\to-\infty,$ тогда как (\ref{4.1}) и (\ref{4.2}) дают
$$
\frac{\theta(z)}{S(z)}=(-1)^{p-s}\frac{C_0}\alpha\exp((C_1-C_2)z)F(z)\prod_{n=s^1+1}^\infty\frac{z^0_n}{z_n}
\prod_{n=p-N+1}^{s^1}\frac1{z_n} \prod_{n=s-N+1}^{s^1}z^0_n,
$$
где $s^1=\max\{s,p\}-N$ и
$$
C_2=\gamma+\sum_{n=s-N+1}^{s^1}\frac1{z^0_n} -\sum_{n=p-N+1}^{s^1}\frac1{z_n} +\sum_{n=s^1+1}^\infty\Big(\frac1{z^0_n} -\frac1{z_n}\Big),
$$
$$
F(z)=z^{p-s}\prod_{n=s-N+1}^{s^1}\frac1{z-z^0_n}\prod_{n=p-N+1}^{s^1}(z-z_n) \prod_{n=s^1+1}^\infty\frac{z-z_n}{z-z^0_n} \to1, \quad
z^2\to-\infty.
$$
При этом для определенности считаем, что
$$
z_{1-N}=\ldots=z_{p-N}=0, \quad z^0_{1-N}=\ldots=z^0_{s-N}=0,
$$
$$
\sum_{n=n_1}^{n_2}(\ldots)=0, \quad \prod_{n=n_1}^{n_2}(\ldots)=1, \quad n_1>n_2.
$$
Таким образом, приходим к соотношениям $C_1=C_2$ и
$$
C_0=(-1)^{p+s}\alpha\prod_{n=s^1+1}^\infty\frac{z_n}{z^0_n} \prod_{n=p-N+1}^{s^1}z_n \prod_{n=s-N+1}^{s^1}\frac1{z^0_n}.
$$
Подставляя их в (\ref{4.1}), получаем представление
$$
\theta(z)=\alpha\exp(\gamma z)\prod_{n=1-N}^{s-N} (z-z_n) \prod_{n=s-N+1}^\infty\frac{z_n-z}{z^0_n} \exp\Big(\frac{z}{z^0_n}\Big),
$$
которое с учетом (\ref{2.3-1}) дает (\ref{2.6}). Теорема 5 доказана.

\medskip
Докажем теорему~6. Согласно лемме~2 соотношения (\ref{4.2-1}) фактически задают, и притом однозначно, некоторую функцию $\theta(z)$ вида
(\ref{2.1}), для которой числа $\{z_n\}_{n\ge1}$ являются нулями с учетом кратности. Помимо них, в соответствии с теоремой~4, эта функция
$\theta(z)$ имеет еще в точности $N$ нулей $\{z_n\}_{n=\overline{1-N,0}},$ которые тоже определяются однозначно. С другой стороны, согласно
уже доказанной теореме~5, полученная функция $\theta(z)$ имеет представление (\ref{2.6}), что и дает первое утверждение теоремы.

Пусть теперь задана полная последовательность $\{z_n\}_{n\ge1-N}$ вида (\ref{2.4}). Обозначим через $\tilde\theta(z)$ функцию, построенную
в первой части доказательства с использованием лишь $\{z_n\}_{n\ge1},$ а через $\{\tilde z_n\}_{n=\overline{1-N,0}}$ -- остальные
получившиеся нули $\tilde\theta(z)$ и рассмотрим функцию
\begin{equation}\label{4.4}
\theta(z):=\tilde\theta(z)\frac{Q_N(z)}{\tilde Q_N(z)}=\frac{Q_N(z)P_N(z)S_0(z)}{\tilde Q_N(z)} +\frac{Q_N(z)}{\tilde Q_N(z)}\int_{-b}^b
w(x)\exp(izx)\,dx,
\end{equation}
где
$$
Q_N(z)=\prod_{n=1-N}^0(z-z_n), \quad \tilde Q_N(z)=\prod_{n=1-N}^0(z-\tilde z_n).
$$
Также рассмотрим единственный полином $Q_{N-1}(z)$ степени меньшей $N,$ такой что нули $\tilde Q_N(z)$ являются нулями и разности
$Q_N(z)P_N(z)-Q_{N-1}(z).$ Очевидно, что степень полинома
\begin{equation}\label{4.5}
P_{N-1}(z):=\frac{Q_N(z)P_N(z)-Q_{N-1}(z)}{\tilde Q_N(z)}-P_N(z),
\end{equation}
в свою очередь, также не превосходит $N-1,$ а функция
\begin{equation}\label{4.6}
f(z):=\frac1{\tilde Q_N(z)}\Big({Q_{N-1}(z)S_0(z)+Q_N(z)}\int_{-b}^b w(x)\exp(izx)\,dx\Big)
\end{equation}
является целой и имеет оценку $f(z)=o(\exp(|{\rm Im\,}z|b))$ при $z\to\infty.$ Кроме того, очевидно, что $f(z)\in L_2(-\infty,\infty).$
Таким образом, в силу теоремы Пэли--Винера (см., например, \cite{5}), она имеет вид
\begin{equation}\label{4.7}
f(z)=\int_{-b}^b \tilde w(x)\exp(izx)\,dx, \quad \tilde w(x)\in L_2(-b,b).
\end{equation}
Сопоставляя формулы (\ref{4.4})--(\ref{4.7}), приходим к (\ref{2.6-1}). Теорема~6 доказана.
\\

{\large\bf Доказательство теоремы 7}
\\

Левое неравенство в (\ref{3.3}) дает оценку
\begin{equation}\label{5.0}
\|\hat w\|_{L_2(-b,b)}^2\le \frac1m\sum_{n=1}^\infty|\hat\theta(z^1_n)|^2.
\end{equation}
Напомним, что последовательность $\{z^1_n\}_{n\ge1}$ является простой и отличается от последовательности $\{z^0_n\}_{n\ge1}$ нулей функции
$S_0(z)$ лишь конечным числом элементов. Без ущерба для общности будем считать, что $z^1_n\ne z^0_n$ тогда и только тогда, когда $n<h$ для
некоторого $h\in{\mathbb N},$ а также что $|\mu_k^N(z^0_k-z^1_n)|\ge5r,$ $n=\overline{1,h-1},$ $k\ge1-N.$

Преобразуем представление (\ref{4.2}) с учетом (\ref{2.3-1}) к виду
$$
S(z)=\alpha \exp(\beta z)\prod_{n=1-N}^\infty\frac{z^0_n-z}{\mu_n}\exp\Big(\frac{z}{\mu_n}\Big),
$$
который вместе с (\ref{2.6}) приводит к соотношению
\begin{equation}\label{4.4-1}
\theta(z)=S(z)\prod_{n=1-N}^\infty\frac{z_n-z}{z^0_n-z}.
\end{equation}
Таким образом, получаем представление
$$
\hat\theta(z^1_n)=S(z^1_n)(d_n-\tilde d_n) =d_nS(z^1_n)\Big(1-\frac{\tilde d_n}{d_n}\Big), \quad n=\overline{1,h-1},
$$
где использованы обозначения
$$
d_n:=\prod_{k=1-N}^\infty\frac{z_k-z^1_n}{z^0_k-z^1_n}, \quad \tilde d_n:=\prod_{k=1-N}^\infty\frac{\tilde z_k-z^1_n}{z^0_k-z^1_n}.
$$
Согласно (\ref{2.4}) условие теоремы дает $|\varkappa_n|\le r$ и $|\tilde\varkappa_n|\le r$ при $n\ge1-N.$ Также имеем
$$
d_n=\prod_{k=1-N}^\infty \Big(1+\frac{\varkappa_k}{\mu_k^N(z^0_k-z^1_n)}\Big), \quad \frac{\tilde d_n}{d_n}
=\prod_{k=1-N}^\infty\Big(1-\frac{\hat\varkappa_k}{\mu_k^N(z_k-z^1_n)}\Big).
$$
Для краткости обозначим $\gamma_{k,n}:=\mu_k^N(z^0_k-z^1_n).$ Тогда $\mu_k^N(z_k-z^1_n)=\gamma_{k,n}+\varkappa_k,$ и наше требование для
$z^1_n$ приобретает вид $|\gamma_{k,n}|\ge5r,$ $n=\overline{1,h-1},$ $k\ge1-N,$ а значит,
$$
\Big|\frac{\varkappa_k}{\mu_k^N(z^0_k-z^1_n)}\Big|= \frac{|\varkappa_k|}{|\gamma_{k,n}|}\le\frac{r}{5r}=\frac15,
$$
$$
\Big|\frac{\hat\varkappa_k}{\mu_k^N(z_k-z^1_n)}\Big|\le
\frac{|\varkappa_k|+|\tilde\varkappa_k|}{|\gamma_{k,n}|-|\varkappa_k|}\le\frac{r+r}{5r-r}=\frac12,
$$
что позволяет записать
$$
d_n=\exp\Big(\sum_{k=1-N}^\infty \ln\Big(1+\frac{\varkappa_k}{\gamma_{k,n}}\Big)\Big), \quad \frac{\tilde
d_n}{d_n}=\exp\Big(\sum_{k=1-N}^\infty \ln\Big(1-\frac{\hat\varkappa_k}{\gamma_{k,n}+\varkappa_k}\Big)\Big).
$$
В результате приходим к оценкам
$$
|d_n| \le \exp\Big(\frac54\sum_{k=1-N}^\infty \Big|\frac{\varkappa_k}{\gamma_{k,n}}\Big|\Big) \le C_r, \quad \Big|1-\frac{\tilde
d_n}{d_n}\Big| \le \sum_{\nu=1}^\infty \frac{2^\nu}{\nu!}\Big(\sum_{k=1-N}^\infty
\Big|\frac{\hat\varkappa_k}{\gamma_{k,n}+\varkappa_k}\Big|\Big)^\nu.
$$
Преобразуем вторую оценку к виду
$$
\Big|1-\frac{\tilde d_n}{d_n}\Big| \le 2\sigma_ne^{2\sigma_n}, \quad \sigma_n:=\sum_{k=1-N}^\infty
\Big|\frac{\hat\varkappa_k}{\gamma_{k,n}+\varkappa_k}\Big|.
$$
Неравенство Коши--Буняковского дает
\begin{equation}\label{5.9}
\left.\begin{array}{l}
 \displaystyle \sigma_n^2\le\alpha_n\|\{\hat\varkappa_k\}\|_{l_2}^2,\\[3mm]
 \displaystyle \alpha_n:=\sum_{k=1-N}^\infty \frac1{|\gamma_{k,n}+\varkappa_k|^2} \le\frac{25}{16}\gamma_n,\\[3mm]
 \displaystyle \gamma_n:=\sum_{k=1-N}^\infty\frac1{|\gamma_{k,n}|^2}\le C,
\end{array}\right\}
\end{equation}
где второе неравенство следует из оценки
$$
\frac1{|\gamma_{k,n}+\varkappa_k|} \le\frac1{|\gamma_{k,n}|}\Big(1-\Big|\frac{\varkappa_k}{\gamma_{k,n}}\Big|\Big)^{-1}
\le\frac1{|\gamma_{k,n}|}\Big(1-\frac{r}{5r}\Big)^{-1}=\frac5{4|\gamma_{k,n}|}.
$$
Итак, приходим к оценкам
\begin{equation}\label{5.9-2}
|\hat\theta(z^1_n)|\le C_r\|\{\hat\varkappa_k\}\|_{l_2}, \quad n=\overline{1,h-1}.
\end{equation}

Пусть теперь $n\ge h.$ Тогда, в силу нашего соглашения, имеем $z^1_n=z^0_n,$ причем $S'(z^0_n)\ne0,$ и справедливо (\ref{2.3}). Подставляя
$z=z^0_n$ в (\ref{4.4-1}), с учетом (\ref{2.4}) получаем
$$
\theta(z^0_n)=-\frac{S'(z^0_n)}{\mu_n^N}d_n\varkappa_n, \quad n\ge h,
$$
где
$$
d_n:=\prod_{{k\ne n}\atop{k=1-N}}^\infty\frac{z_k-z^0_n}{z^0_k-z^0_n} =\prod_{{k\ne
n}\atop{k=1-N}}^\infty\Big(1+\frac{\varkappa_k}{\mu_k^N(z^0_k-z^0_n)}\Big), \quad n\ge h.
$$
Таким образом, приходим к соотношению
\begin{equation}\label{5.1}
\hat\theta(z^0_n)=-\frac{S'(z^0_n)}{\mu_n^N}(d_n\hat\varkappa_n+\hat d_n\tilde\varkappa_n), \quad n\ge h.
\end{equation}
Положив для краткости $\gamma_{k,n}:=\mu_k^N(z^0_k-z^0_n),$ факторизуем $d_n$ следующим образом:
\begin{equation}\label{4.6-1}
d_n=d_{n,1}d_{n,2}, \quad d_{n,1}:=\prod_{k\in\Omega_r(n)}\Big(1+\frac{\varkappa_k}{\gamma_{k,n}}\Big), \quad
d_{n,2}:=\prod_{k\in\Theta_r(n)}\Big(1+\frac{\varkappa_k}{\gamma_{k,n}}\Big),
\end{equation}
где
$$
\Omega_r(n):=\{k:k+N\in{\mathbb N},\,0<|\gamma_{k,n}|<5r\},
$$
$$
\Theta_r(n):=\{k:k+N\in{\mathbb N},\,|\gamma_{k,n}|\ge5r\}.
$$
Имеем
\begin{equation}\label{4.8}
|d_{n,1}|\le\Big(1+\frac{r}\chi\Big)^R, \quad \chi:=\min_{k\ne n}|\gamma_{k,n}|,  \;\; R:=\max_{n\ge h}\#\Omega_r(n)<\infty,
\end{equation}
$$
d_{n,2}=\exp\Big(\sum_{k\in\Theta_r(n)}\ln\Big(1+\frac{\varkappa_k}{\gamma_{k,n}}\Big)\Big)
$$
и, следовательно,
\begin{equation}\label{4.9}
|d_{n,2}|\le\exp\Big(\frac54\sum_{k\in\Theta_r(n)}\Big|\frac{\varkappa_k}{\gamma_{k,n}}\Big|\Big)\Big) < \exp({2r\sqrt{\gamma_n}}), \quad
\gamma_n:=\sum_{k\in\Theta_r(n)}\frac1{|\gamma_{k,n}|^2}\le C,
\end{equation}
где $C$ не зависит от $n.$ Итак, получаем $|d_n|\le C_r$ для всех $n\ge h.$

Займемся вторым слагаемым в скобках в (\ref{5.1}). Имеем
\begin{equation}\label{5.7}
\hat d_n=d_{n,2}\Big(\hat d_{n,1}+ \tilde d_{n,1}\Big(1-\frac{\tilde d_{n,2}}{d_{n,2}}\Big)\Big),
\end{equation}
$$
\frac{\tilde d_{n,2}}{d_{n,2}}=\prod_{k\in\Theta_r(n)}\frac{\gamma_{k,n}+\tilde\varkappa_k}{\gamma_{k,n}+\varkappa_k}
=\prod_{k\in\Theta_r(n)}\Big(1-\frac{\hat\varkappa_k}{\gamma_{k,n}+\varkappa_k}\Big),
$$
$$
\Big|\frac{\hat\varkappa_k}{\gamma_{k,n}+\varkappa_k}\Big|\le\frac{2r}{5r-r}=\frac12.
$$
Таким образом, приходим к оценке
$$
\Big|1-\frac{\tilde
d_{n,2}}{d_{n,2}}\Big|=\Big|1-\exp\Big(\sum_{k\in\Theta_r(n)}\ln\Big(1-\frac{\hat\varkappa_k}{\gamma_{k,n}+\varkappa_k}\Big)\Big)\Big|=
$$
$$
=\Big|\sum_{\nu=1}^\infty\frac1{\nu!}
\Big(\sum_{k\in\Theta_r(n)}\ln\Big(1-\frac{\hat\varkappa_k}{\gamma_{k,n}+\varkappa_k}\Big)\Big)^\nu\Big|
\le\sum_{\nu=1}^\infty\frac{2^\nu}{\nu!}\Big(\sum_{k\in\Theta_r(n)} \Big|\frac{\hat\varkappa_k}{\gamma_{k,n}+\varkappa_k}\Big|\Big)^\nu=
$$
\begin{equation}\label{5.8}
=\sum_{\nu=1}^\infty \frac{(2\sigma_n)^\nu}{\nu!}\le 2\sigma_n\exp(2\sigma_n), \quad \sigma_n:=\sum_{k\in\Theta_r(n)}
\Big|\frac{\hat\varkappa_k}{\gamma_{k,n}+\varkappa_k}\Big|.
\end{equation}
Аналогично (\ref{5.9}) получаем неравенства
\begin{equation}\label{5.9-1}
\sigma_n^2\le\alpha_n\|\{\hat\varkappa_k\}\|_{l_2}^2, \quad \alpha_n:=\sum_{k\in\Theta_r(n)} \frac1{|\gamma_{k,n}+\varkappa_k|^2}
\le\frac{25}{16}\gamma_n^2,
\end{equation}
где $\gamma_n$ определено в (\ref{4.9}). Оценим теперь первое слагаемое в скобках в (\ref{5.7}):
$$
\hat d_{n,1}=\sum_{k\in\Omega_r(n)}\frac{\hat\varkappa_k}{\gamma_{k,n}}
\prod_{{\nu\in\Omega_r(n)}\atop{\nu<k}}\Big(1+\frac{\tilde\varkappa_k}{\gamma_{k,n}}\Big)
\prod_{{\nu\in\Omega_r(n)}\atop{\nu>k}}\Big(1+\frac{\varkappa_k}{\gamma_{k,n}}\Big).
$$
По аналогии с (\ref{4.8}) имеем
\begin{equation}\label{5.10}
|\hat d_{n,1}|\le\frac1{\chi}\Big(1+\frac{r}\chi\Big)^{R-1}\sum_{k\in\Omega_r(n)}|\hat\varkappa_k|.
\end{equation}
Наконец, в силу (\ref{2.1}) и (\ref{2.3}) справедлива оценка $S'(z^0_n)=O(\mu_n^N),$ $n\to\infty.$

Итак, используя (\ref{5.1})--(\ref{5.10}), приходим к неравенству
$$
|\hat\theta(z^0_n)|\le C_r\Big(|\hat\varkappa_n|+\|\{\hat\varkappa_k\}\|_{l_2}|\tilde\varkappa_n|
+\sum_{k\in\Omega_r(n)}|\hat\varkappa_k|\Big), \quad n\ge h,
$$
откуда получаем оценку
\begin{equation}\label{5.11}
\sqrt{\sum_{n=h}^\infty|\hat\theta(z^0_n)|^2}\le C_r\Big(\|\{\hat\varkappa_k\}\|_{l_2}
+\sqrt{\sum_{n=h}^\infty\Big(\sum_{k\in\Omega_r(n)}|\hat\varkappa_k|\Big)^2}\Big).
\end{equation}
Поскольку
$$
\sum_{n=h}^\infty\Big(\sum_{k\in\Omega_r(n)}|\hat\varkappa_k|\Big)^2 \le R\sum_{n=h}^\infty \sum_{k\in\Omega_r(n)}|\hat\varkappa_k|^2 =
R\sum_{k=1-N}^\infty|\hat\varkappa_k|^2\sum_{n\in\Omega_r^*(k)}1,
$$
где $R$ определено в (\ref{4.8}), а $\Omega_r^*(k)=\{n:n\ge h,\,k\in\Omega_r(n)\},$ второе слагаемое в правой части (\ref{5.11}) не
превосходит $\sqrt{RR^*}\|\{\hat\varkappa_k\}\|_{l_2},$ где $R^*:=\displaystyle\max_{k\ge1-N}\#\Omega^*_r(k).$ Таким образом, с учетом
соотношений $z^1_n=z^0_n$ при $n\ge h,$ а также $\hat\varkappa_n\equiv\mu_n^N\hat z_n$ при $n\ge1-N,$ неравенства (\ref{5.0}),
(\ref{5.9-2}) и (\ref{5.11}) дают оценку (\ref{2.8}). Теорема~7 доказана.


\end{document}